\documentclass[11pt]{article}
\usepackage{graphics, amsfonts, amsmath, amsthm, amstext, amssymb, amscd, color, rotating, makeidx}
\usepackage{calligra,mathrsfs} 
\usepackage{ulem} 
\usepackage{tikz-cd} 
\usepackage[square,sort,comma,numbers]{natbib}
\usepackage{hyperref}
\hypersetup{
    colorlinks, 
    citecolor=red,
    filecolor=black,
    linktoc=all,     
    linkcolor=blue, 
    urlcolor=blue 
}
\usepackage[margin=0.9in]{geometry} 
\usepackage{fancyhdr}  
\newtheorem{theorem}{Theorem}[section]
\newtheorem{definition}[theorem]{Definition}
\newtheorem{lemma}[theorem]{Lemma}
\newtheorem{proposition}[theorem]{Proposition}

\newtheorem{corollary}[theorem]{Corollary}

\newtheorem{conjecture}[theorem]{Conjecture}
\newtheorem*{remark}{Remark}

\newcommand{\barD}{\overline{\partial}}
\newcommand{\ddbar}{\partial\overline{\partial}} 
\newcommand{\supp}{\operatorname{supp}} 
\newcommand{\Hol}{\operatorname{Hol}} 
\newcommand{\Jac}{\operatorname{Jac}} 
\newcommand{\Ric}{\operatorname{Ric}} 
\newenvironment{claim}[1][]{\par\smallbreak
   \noindent \underline{\textsc{Claim}}. \rmfamily}{\smallbreak} 
\newenvironment{cproof}[1][Proof]{%
  \proof[#1]%
  }{\endproof}

\usepackage[]{amsmath, amsthm, amsfonts, verbatim, amssymb,pb-diagram, stmaryrd, wasysym,tikz}
\usepackage{epsfig}

\def\cF{\mathcal F}

\def\cO{\mathcal O}

\def\cT{{\mathcal T}}

\def\bC{{\mathbb C}}

\def\cT{{\mathcal T}}
\def\cF{{\mathcal F}}

\def\cU{\mathcal U}

\def\tX{\widetilde X}
\def\tY{\widetilde Y}

\def\pd{\partial}

\def\oz1{d{\overline z}^1}
\def\oz2{d{\overline z}^2}
\def\oz3{d{\overline z}^3}

\def\oV{\overline V}
\def\oW{\overline W}
\def\oI{\overline I}

\def\oz{\overline z}

\def\op{\overline D}
\def\oIq1{\oI_1\cdots\oI_{q-1}}
\def\oIq2{\oI_1\cdots\oI_{q-2}}
\def\op{\overline \partial}

\def\tV{\widetilde{V}}

\def\oZ{\overline Z}

\def\ddz{\frac{\pd}{\pd z}}

\def\ddw{\frac{\pd}{\pd w}}
\def\ddz1{\frac{\pd}{\pd z^1}}
\def\ddw1{\frac{\pd}{\pd w^1}}
\def\oddw1{\overline{\frac{\pd}{\pd w_1}}}

\def\tZ{\widetilde{Z}}
\def\oZ{\overline Z}
\def\oX{\overline X}

\def\tk{\widetilde k}
\def\tU{\widetilde U}
\def\tV{\widetilde V}
\def\tcU{\widetilde {\mathcal U}}

\def\ty{\widetilde y}
\def\tS{\widetilde S}
\def\hs{\widehat s}

\def\toZ{\widetilde{\oZ}}

\begin{document}

\title{ Quasi-projective Manifolds Uniformized by Carath\'eodory Hyperbolic Manifolds and Hyperbolicity of Their Subvarieties}

\author{Kwok-Kin Wong and Sai-Kee Yeung
\thanks{The second author was partially supported by a grant from the National Science Foundation}}
\maketitle
\begin{abstract}  
{\it   Let $M$ be a Carath\'eodory hyperbolic complex manifold. We show that $M$ supports a  real-analytic bounded strictly plurisubharmonic function. If $M$ is also complete K\"ahler, we show that $M$ admits the Bergman metric. When $M$ 
is strongly Carath\'eodory hyperbolic and 
is the universal covering of a quasi-projective manifold $X$, the Bergman metric  can be estimated in terms of a Poincar\'e type metric on $X$. It is also
proved that any quasi-projective (resp. projective) subvariety of $X$ is of log-general type (resp. general type),
a result consistent with a conjecture of Lang. 
}
\end{abstract}

\section*{Introduction}

Invariant pseudo-metrics are important tools for studying complex manifolds. When the pseudo-metrics are nondegenerate, i.e. they actually define metrics, the underlying complex  manifolds usually possess very interesting properties. 
As a generalization of the Poincar\'e metric on the unit disk $\Delta\subset \mathbb{C}$, for any complex manifold $X$, one can define the infinitesimal Kobayashi-Royden pseudo-metric
\[
\mathcal{K}_X(v):= \inf\{ \lambda>0 \mid \exists f:\Delta\rightarrow X, f(0)=p, \lambda f'(0)=v\}, \quad v\in T_pX,
\]
which is integrated to a pseudo-distance
\[
d_{K,X}(p,q):= \inf \{ \int_\gamma \mathcal{K}_X(\gamma'(t))dt \mid  \gamma \text{ is a piecewise smooth courve joining $p$ to $q$} \}.
\]
A complex manifold $X$ is said to be \textit{Kobayashi hyperbolic} if $d_{K,X}$ is nondegenerate. 

Another important invariant pseudo-metric is the Bergman pseudo-metric $\mathcal{B}_X$, which is constructed by considering square-integrable holomorphic top forms (see \S \ref{Berg} for a brief review of the construction in the K\"ahler case). In this article, $\mathcal{B}_X$ will be closely related to  our main focuses, the \textit{infinitesimal Carath\'eodory pseudo-metric} $E_X$  and the \textit{Carath\'eodory pesudo-distance} $d_{C,X}$. They are in some sense dual to $\mathcal{K}_X$ and $d_{K,X}$ respectively. 
For any bounded domain $\Omega\subset \mathbb{C}^n$, the pseudo-metrics $ \mathcal{K}_\Omega, \mathcal{B}_\Omega,E_\Omega$ are always nondegenerate. 
For arbitrary complex manifolds, they are only pseudo-metrics and many properties among the three metrics known for bounded domains are not clear for complex manifolds. 
For the constructions and basic properties of these metrics, see for example \cite{Koba1998}.

We say that a complex manifold $X$ is\textit{ Carath\'eodory hyperbolic} if the infinitesimal Carath\'eodory pseudo-metric $E_X$ is nondegenerate. Following from definition, such manifolds are noncompact (see \S \ref{caradis} for definition and its related notions). 
Our first result is

\begin{theorem}
\label{psh}
Any Carath\'edory hyperbolic complex manifold
admits a bounded  real-analytic strictly plurisubharmonic function.
\end{theorem}
\noindent
A result of similar nature may be known among experts, but since we cannot locate a precise reference and the result is to be used later in the paper, we provide a detailed proof here.
If furthermore the Carath\'eodory pseudo-distance is nondegenerate and complete, which we named strongly Carath\'eodory hyperbolic (Definition \ref{scara}), we obtain from Theorem \ref{psh} that the manifold is Stein (Corollary \ref{exhaust}). 

Now suppose the Carath\'eodory hyperbolic complex manifold can be given a complete K\"ahler metric, 
we can deduce from Theorem \ref{psh} and the method of $L^2$-estimates that
\begin{theorem}
For any Carath\'eodory hyperbolic complete K\"ahler manifold, the Bergman pseudo-metric is nondegenerate. 
\label{bergm}
\end{theorem}
\noindent 
In fact, once we have the existence of a smooth bounded plurisubharmonic exhaustion function, Theorem \ref{bergm} follows from $L^2$-estimates as given in the proof of  \cite[Corollary  1 (a)]{Yeun2003}, \cite[Theorem 1]{Yeun2000} and \cite[Proposition 2.3]{Chen2003}. To be more self-contained, we will add the details (Proposition \ref{genM}) in this article.
Theorem \ref{bergm} allows us to remove a nondegeneracy assumption in a classical comparison theorem between infinitesimal Carath\'eodory metric and Bergman metric for Carath\'eodory hyperbolic complete K\"ahler manifolds (Corollary \ref{ineqcb}). We expect this to be useful for studying some unbounded domains. 

\bigskip

We continue to look at projective and, more generally, quasi-projective manifolds whose universal covers are Carath\'eodory hyperbolic. Our main motivation is the following  well-known question from complex hyperbolicity:
\begin{conjecture}[Lang \cite{Lang1986}]
A projective variety $X$ is Kobayashi hyperbolic if and only if all subvarieties of $X$ (including $X$ itself) are of general type.
\label{langconj}
\end{conjecture}
Lang's Conjecture \ref{langconj} speculates that there are intimate relationships between the analytic notion of Kobayashi metric and the complex differential and algebro-geometric notion of general type. Although Conjecture \ref{langconj} is largely open, it is believed that the necessity condition for all subvarieties to be of general type is guaranteed by various analytic or differential geometric hyperbolicity conditions. 
It is well-known that Carath\'eodory hyperbolic implies Kobayashi hyperbolic.
In this spirit, we are going to consider Carath\'eodory metric instead of Kobayashi metric.  We first establish
\begin{theorem}
Let $X$ be a compact K\"ahler manifold whose universal cover is Carath\'edory hyperbolic. Then all subvarieties of $X$ (including $X$ itself) are projective varieties of general type.
\label{gtcpt}
\end{theorem}

In fact, our main interest here is how Lang's Conjecture \ref{langconj} could be generalized to quasi-projective manifolds. As is already pointed out in \cite{Lang1986}, it is natural to consider various notions of hyperbolicity modulo a subset, i.e. those conditions are satisfied outside certain exceptional subsets.
For quasi-projective manifolds, a first candidate for such subset would be the compactifying divisor of a given projective compactification. This points to the study of variety of log-general types.
We refer to, for example, the monograph \cite{Vojta1987}, the recent surveys \cite{Java2020, AT2020} 
and references therein for some account of this rich subject, especially the connections to arithmetic. 
 For quasi-projective manifolds, an immediate modification should be the replacement of the notion of general type by log-general type. While subvarieties of projective manifolds are automatically projective, there are examples of quasi-projective manifolds that contain analytic subvarieties which are not quasi-projective (Proposition \ref{egnonqp}). Therefore, we propose the following analogue of Conjecture \ref{langconj}:
\begin{conjecture}
A quasi-projective variety $X$ is Kobayashi hyperbolic if and only if all quasi-projective subvarieties of $X$ are of log-general type and all projective subvarieties are of general type.
\label{conj}
\end{conjecture}
Similar to the compact case, we are going to deduce our main result
\begin{theorem}[Main Theorem]
Let $X$ be a  quasi-projective manifold whose universal covering is strongly Carath\'eodory hyperbolic. Then all quasi-projective subvarieties of $X$ (including $X$ itself) are of log-general type
and all projective subvarieties are of general type.
\label{gtncpt}
\end{theorem}

\noindent 
As in the proof of Theorem \ref{gtcpt}, we consider the Bergman metric on the universal covering obtained from Theorem \ref{bergm}, which descends to a Hermitian metric on $K_X$ with constant negative Ricci curvature. Crucial to us is the observation that this Hermitian metric can be estimated and compared with a Poincar\'e type metric that is naturally defined for any quasi-projective manifold. 
We will show that such estimate can be obtained if the universal covering supports a bounded plurisubharmonic exhaustion function (Proposition \ref{ker-est}), which is a consequence of being strongly Carath\'eodory hyperbolic (Proposition \ref{exhaust}). Such estimates are also of interest in their own right. 
\bigskip

 We briefly outline some major techniques in this article.
The machinery includes the method of $L^2$-estimates of $\barD$, the construction of appropriate weight functions and the utilization of the Bergman metrics. There are two main ingredients.  The first is the construction of a bounded {\it strictly} plurisubharmonic function (Theorem \ref{psh}).
The second is the estimates and construction of an appropriate Hermitian metric for the canonical line bundle on the regular part of a possibly singular quasi-projective subvariety 
of the ambient manifold so that enough log-canonical sections on a compactification can be guaranteed.

The existence of a bounded strictly plurisubharmonic function 
proved here depends on the properties  of  the Carath\'edory metric and distances.  
Earlier some geometric consequences of the existence of a bounded plurisubharmonic function have been given
in  \cite[Theorem 4]{Yeun2000} and \cite[Corollary 1]{Yeun2003}.  The paper here gives further geometric applications.  

When combined with the method of $L^2$-estimates, it is shown that the $L^2$-holomorphic sections of the canonical bundle $K_M$ of a Carath\'edory hyperbolic complete K\"ahler manifold $M$ separate points and generate jets of $K_M$ (Proposition \ref{genM}). The same strategy can be applied to pluricanonical bundles $qK_M$ for any $q>0$ (Lemma \ref{genqM}). 
Now suppose $X$ is a compact K\"ahler manifold with Carath\'eodory hyperbolic universal cover $M$. 
Then the Bergman metric on $M$ is well-defined and descends to $X$, forming a complete K\"ahler metric of strictly negative Ricci curvature. By Kodaira's embedding theorem, $K_X$ is ample. Thus $X$ is of general type. For a subvariety $Z\subset X$, one needs to show that a resolution $\widetilde{Z}\rightarrow Z$ is of general type. 
Let $Y$ be an irreducible component of the pull-back of $Z$ to $M$.
Since $M$ is Carath\'eodory hyperbolic, it has a bounded strictly plurisubharmonic function $\varphi$.
One can pull-back the restriction $\varphi|_{Y}$
from $Y$ to its resolution $\widetilde{Y}\rightarrow Y$, but it is strictly plurisubharmonic only outside the inverse image of the singular set of $\tY$. Thus one cannot simply repeat the argument for $X$.  For this, we supply an argument using Poincar\'e series to show that $\tZ$ must be of general type.   Alternatively, the reader will observe that the following more elaborate argument for the quasi-projective setting works for the projective setting as well.

When $X$ is a quasi-projective manifold, to prove that $X$ is of log-general type, one needs to show that there is a compactification $\overline{X}\supset X$ such that $K_{\overline{X}}+ [D]$ is big, where $D=\overline{X}-X$ can be chosen to be a divisor of simple normal crossing. To do this, we will use a technique of Mok \cite{Mok1986}, which requires 
positivity and a uniform bound of the curvature of $K_X$ with respect to some Hermitian metric. We construct such metric out of the Bergman kernel on its universal cover $M$ and the natural metric of Poincar\'e type on $X$. Here the bounded geometry of a Poincar\'e type metric will play an essential role. For quasi-projective subvarieties $Z\subset X$, we have similar difficulties as in the compact case that the bounded strictly plurisubharmonic on $M$ does not pull-back to a strictly plurisubharmonic function on all of the cover $\tY\rightarrow \tZ$ for the resolution $\tZ\rightarrow Z$. The key observation is that although we only get some $\psi$ that is strictly plurisubharmonic on $\tY-B$ for some subset $B\subset \tY$, the estimates of the Poincar\'e metric by the Bergman metric constructed from $\psi$ 
(Proposition \ref{ker-est}) 
will allow us to obtain the desired curvature bound on $K_{\tZ-A}$ for some subset $A\subset \tZ$. We show that this is sufficient for us to apply Mok's technique \cite{Mok1986}  to conclude that the space of $L^2$-holomorphic sections of $K_{\tZ-A}$ is of the desired growth. By combining these with an estimate near the boundary, we show that those $L^2$ sections extend to boundary, hence concluding Theorem \ref{gtncpt}. 
Note that we need the exhaustion coming from strongly Carath\'eodory hyperbolicity to obtain the estimates of Bergman metric with respect to the pull-back Poincar\'e type metric on universal cover.
\bigskip

We mention some previous works closely related to us. They focus exclusively on compact manifold quotients of the form $X=M/\Gamma$. In \cite{Wu1987}, $M$ is a Carath\'eodory hyperbolic complex manifold. It is shown that one can construct an invariant volume form $v_{C,M}$ on $M$, which descends to $X$ to give a smooth Hermitian metric with negative Ricci curvature on $X$. The volume form $v_{C,M}$ is in general only a pseudo-volume form. Nowadays, a complex manifold with nondegenerate $v_{C,M}$ is usually said to be Carath\'eodory measure hyperbolic, cf. \cite{Eise1970}. In \cite{Kiku2011}, assuming $M$ is Carath\'eodory measure hyperbolic, it is shown that there is a smooth bounded plurisubharmonic function on $M$ which is strictly plurisubharmonic outside the set of degeneracy of $v_{C,M}$. Combining with \cite[Proposition 2.3]{Chen2003} (which is the same as Theorem \ref{bergm}) and \cite{Bou2002}, this implies the positivity and some explicit volume estimates of $K_X$ in terms of the volume form on $X$ obtained by descending $v_{C,M}$. Similar estimates are shown to hold in \cite{Ng2018} for all subvarieties of $X$ in the setting of \cite{Kiku2011}, thus solving the conjecture of Kikuta in \cite{Kiku2013}. Note also that these volume estimates implies the bigness of the corresponding canonical bundles, hence showing that all subvarieties of the compact manifold $X$ are of general type. The compactness is essential in the techniques of Wu \cite{Wu1987}, Kikuta \cite{Kiku2011,Kiku2013} and Ng \cite{Ng2018}. For quasi-projective case, we believe that an analysis in analogues to Proposition \ref{ker-est} near the boundary is necessary. 

%
%

\tableofcontents

\section{Carath\'eodory  hyperbolicity and geometric consequences}

\subsection{Poincar\'e distance on unit disk}
\label{Poin}
Let $\Delta$ be the unit disk on complex plane. 
For $z\in \Delta$ and a holomorphic vector $v\in T_z\Delta$, the (infinitesimal) Poincar\'e metric is defined by
\[
ds^2_\Delta = \frac{dz\otimes d\overline{z}}{(1-|z|^2)^2}
\]
with the induced norm
\[
\|v\|_P := \frac{|v|}{1-|z|^2},
\]
where $v$ is identified to a complex number via $T_z\Delta\cong \mathbb{C}$. Note that  $\|v\|_P>0$ if and only if $|v|>0$. By integrating over the line segment joining $0$ and $z$ with respect to $ds^2_\Delta$, we obtain
the Poincar\'e distance $\ell_P(0,z)$ between $0$ and $z\in \Delta$, which is given by
\[
\ell_P(z):=\ell_P(0,z)=\frac{1}{2}\log \frac{1+|z|}{1-|z|}=\frac{1}{2}\log \frac{(1+|z|)^2}{1-|z|^2}.
\]

\subsection{Carath\'eodory metric and distance}
\label{caradis}

Let $M$ be a complex manifold. Denote by $\Hol(M,\Delta)$ the set of holomorphic maps from $M$ to $\Delta$.
The \textit{infinitesimal Carath\'eodory  pesudo-metric} on $M$ is defined by
\begin{equation}
E_M(x,v):=\sup\{ \|df_x(v)\|_P \mid f\in \Hol(M,\Delta), f(x)=0\}, \quad x\in M, \quad v\in T_xM.
\label{em}
\end{equation}
Note that $E_M\geq 0$ in general.
If the set of bounded nonconstant holomorphic functions on $M$ is nonempty, then from normal family argument (cf. \cite[Chapter 4.2]{Koba1998}), $E_M(x,v)$ is realized by some holomorphic mapping $f_{x.v}\in \Hol(M,\Delta)$ for any given $x\in M$ and $v\in T_xM$. 

\begin{definition}
Let $M$ be a complex manifold.
We say that $E_M$ is nondegenerate at $x\in M$ if $E_M(x,v)>0$ for every $v\in T_xM-\{0\}$. 
$M$ is said to be {\it Carath\'eodory hyperbolic} if $E_M$ is nondegenerate for all $x\in M$. In this case, $E_M$ is simply called the infinitesimal Carath\'eodory metric of $M$.
\label{cara}
\end{definition} 
\begin{remark}
\begin{enumerate}
\item By definition, if $M$ has no nonconstant bounded holomorphic functions, then $M$ is not Carath\'eodory hyperbolic. Thus a Carath\'eodory hyperbolic complex manifold is necessarily noncompact.
\item We will be interested in complex manifolds $X$ whose universal covering $\tX$ is Carath\'eodory hyperbolic. As $E_{\tX}$ is invariant under biholomorphism, we may in principle say that a complex manifold $X$ is `Carath\'eodory hyperbolic' if the $E_{\tX}$ is nondegenerate.  However this notion should not be confused with the statement that the infinitesimal Carath\'eodory pseudo-metric $E_X$ on
$X$ is nondegenerate, since there could be no nonconstant bounded holomorphic functions on $X$.
\end{enumerate}
\end{remark}

Our notion of Carath\'eodory hyperbolicity in Definition \ref{cara} may not be the conventional one, such as that in \cite{Koba1998}.
Recall that 
\begin{equation}
d_{C,M}(p,q):= \sup \{ \ell_P(f(p),f(q)) \mid f\in \Hol(M,\Delta)\} ,\quad p,q\in M,
\label{dcm}
\end{equation}
defines the \textit{Carath\'eodory pseudo-distance} on $M$. It is said to be nondegenerate   if $d_{C,M}(p,q)> 0$ whenever $p\neq q$. 
In \cite[p. 174]{Koba1998}, Kobayashi defined for a complex space $M$ to be `Carath\'eodory hyperbolic', or `C-hyperbolic' based on the nondegeneracy of $d_{C,M}$
instead of the infinitesimal Carath\'eodory pseudo-metric $E_M$. 
In general, the nondegeneracy of $d_{C,M}$ and $E_M$ are not necessarily equivalent.

\begin{definition}
\label{scara}
A complex manifold $M$ is said to be strongly Carath\'eodory hyperbolic if it is Carath\'eodory hyperbolic and $d_{C,M}$ is complete nondegenerate.  
\end{definition}

There are a lot of interesting examples of strongly Carath\'eodory hyperbolic manifolds. They include bounded symmetric domains, Teich\"uller space $\cT_{g,n}$ of
hyperbolic Riemann surfaces with genus $g$ and $n$-punctures, and more generally all bounded domains in $\mathbb{C}^n$. Some unbounded domains are also known to be strongly Carath\'eodory hyperbolic, cf. \cite{AGK2016}. On the other hand, we will show that strongly Carath\'eodory hyperbolic implies Stein, see Corollary \ref{exhaust}.


\subsection{Bounded plurisubharmonic function}

\begin{theorem}[= Theorem \ref{psh}]
\label{psh'}
Let $M$ be a Carath\'edory hyperbolic complex manifold. Then
there exists a bounded real-analytic  strictly plurisubharmonic function $\varphi$ on $M$.
\end{theorem}

\begin{proof}
Assume there is a countable family $\{f_j\}_{j\in J}\subset \Hol(M,\Delta)$ such that for each $x\in M$,\\

\begin{tabular}{cl}
$(\sharp)$: & there exist $n=\dim M$ members $\{h_1,\dots,h_n\}\subset \{f_j\}_{j\in J}$ (depending on $x$) so that  \\
			& the Jacobian determinant of $(h_1,\dots,h_n)$ at $x$ is nondegenerate.
\end{tabular}\\

Since $|f_j|<1$ for all $j\in J$, we can find $\{\alpha_{j}\in \mathbb{R}_{>0}: j\in J\}$ such that the sum
\[
\varphi(x):= \sum_{j\in J} \alpha_j^2 |f_{j}(x)|^2<1 
\]
for any $x\in M$. It follows that $\varphi$ converges uniformly on compact subsets of $M$ to a bounded real-analytic function on $M$. Then
\[
\sqrt{-1}\ddbar \varphi =   \sum_{j\in J} \alpha_j^2 \sqrt{-1} \partial f_j\wedge \barD \, \overline{f_j}.
\]
By Cauchy's estimates, $\sqrt{-1}\ddbar \varphi$ converges uniformly on compact subsets of $M$. Now by ($\sharp$), $\varphi$ is strictly plurisubharmonic on $M$. Thus it suffices to construct a countable family $\{f_j\}_{j\in J}\subset \Hol(M,\Delta)$ such that $(\sharp)$ holds for every $x\in M$.

Let $x\in M$ be fixed. Since $M$ is Carath\'eodory hyperbolic, there exists bounded holomorphic functions $f_{1,1},\dots,f_{1,n} \in \Hol(M,\Delta)$ such that the Jacobian determinant $\det \Jac f_x$ of the holomorphic map $f_x:=(f_{11}, \dots  ,f_{1n}):M\rightarrow \mathbb{B}^n$ is nondegenerate at $x$. Thus for any $\beta_1,\dots,\beta_n\in \mathbb{R}_{>0}$, $(\sqrt{-1}\ddbar \sum_{j=1}^n \beta_j^2 |f_{1j}|^2)(x)>0$.
 
 Let 
 \[
 V_x=\{y\in M: (\det \Jac f_x)(y))=0\}.
 \]  
Since $V_x\subset M$ is an analytic subvariety, it has at most countably many irreducible components and each of them has dimension $\leq n-1$. Name the  $f_x$ constructed above as $f_1$.
 Let $V_0\subset V_x$ be an irreducible component  and $w\in V_0$ be a fixed smooth point.  Repeat the above argument with $x$ replaced by $w$. We get another holomorphic map $f_2:=f_w: M\rightarrow \mathbb{B}^n$, 
such that there exist $\beta_1,\beta_2>0$ with
$\det \Jac(\beta_1 f_1+\beta_2 f_2)$ is nondegenerate on $(M-V_x)\cup (V_x-V_w)$. Note that the dimension of any irreducible component of $V_x\cap V_w$ is $\leq n-2$. 
 Repeat the above argument for each irreducible component $V_i$ of $V$ ($I\subset \mathbb{N}$)
, we obtain a countable family $\{f_{i}=(f_{i1},\dots,f_{in})\}_{i\in I}\subset \Hol(M,\mathbb{B}^n)$, where $f_{ij}\in \Hol(M,\Delta)$ for all $i\in I$ and $j=1,\dots,n$. 
For the family $\{f_{ij}\}$, we see that $(\sharp)$ holds outside a subvariety $S\subset M$ whose countably many irreducible components are of dimension $\leq n-2$. Now by induction, we can enlarge the family $\{f_{ij}\}$ to a countable family, which we denote without loss of generality by $\{f_j\}_{j\in J}\subset \Hol(M,\Delta)$, such that $(\sharp)$ holds for any $x\in M$.
\end{proof}

Recall that a complex space is Stein if and only if it admits a continuous strictly plurisubharmonic exhaustion function (\cite{Gra1958,Nara1962}). From Theorem \ref{psh}, we have

\begin{corollary} 
Strongly Carath\'eodory hyperbolic complex manifolds are Stein.
\label{stein}
\end{corollary}
\begin{proof}
Let $x_0\in M$ be fixed. Given $x\in M$, recall from definition that $d_{C,M}(x_0,x)=\sup_h\ell_P(h(x))$ among all holomorphic map $h:M\rightarrow \Delta$ such that $h(x_0)=0$. Note that $\ell_P(0)=0$. For $z\in \Delta-\{0\}$, we have $\ell_P(z)>0$ and
\begin{align*}
\sqrt{-1}\pd\op \ell_P(z)
&=\frac{1+|z|^2}{4(1-|z|^2)^2|z|}|dz|^2 \\
&>\bigg|\frac{|z|}{2z(1-|z|^2)}\bigg|^2|dz|^2 
= \sqrt{-1}\partial \ell_P\wedge \barD\ell_P (z),
\end{align*}
which implies that
\begin{align*}
\sqrt{-1}\ddbar \ell_P^2 (z)
&=2\sqrt{-1}\partial \ell_P\wedge \barD\ell_P (z)+2\ell_P(z)\sqrt{-1}\pd\op \ell_P(z) \nonumber \\
&> (2+2\ell_P(z)) \sqrt{-1}\partial \ell_P\wedge \barD\ell_P(z).
\end{align*}
So $\ell_P^2$ satisfies the Sub-Mean-Value Property and thus is subharmonic. It follows that $d_{x_0}(x):=(d_{C,M}(x_0,x))^2=\sup_h\ell_P^2(h(x))$
 is a plurisubharmonic function on $M$.
Now $d_{C,M}$ is complete implies $d_{x_0}$ is an unbounded exhaustion. We can approximate $d_{x_0}$ by smooth ones by standard methods.  Hence there a smooth plurisubharmonic exhaustion $d'_{x_0}$. Let $\varphi$ be given by Theorem \ref{psh'}. Then $\rho:=\varphi+d'_{x_0}$ is as desired. 
\end{proof}

We  recall that a complex manifold is said to be hyperconvex if it admits a plurisubharmonic exhaustion function which is bounded from above. The following result will be used in next section.

\begin{proposition}
Strongly Carath\'eodory hyperbolic complex manifolds are hyperconvex.

\label{exhaust}
\end{proposition}

\begin{proof}
Let $M$ be a strongly Carath\'eodory hyperbolic complex manifold.
Fix a point $p\in M$. Define the family $\mathcal{F}=\{ f\in \Hol(M,\Delta): f(p)=0\}$. 

For each $h\in \mathcal{F}$, $h^*\ell_P$ is plurisubharmonic on $M$ with poles at the zero locus of $h$. Let $a>0$ be sufficiently small so that the disk $\Delta(0,2a)\subset \Delta$. Let $\rho_a$ be a cutoff function on $\Delta$, where $\rho_a\equiv 1$ on $h(M)-\Delta(0,2a)$ and $\equiv 0$ on $\Delta(0,a)$. Then there is $\lambda>0$ such that 
\begin{equation}
\label{hlp2}
h^*( \lambda |z|^2+\rho_a(z)\ell_P(z))
\end{equation}
is strictly plurisubharmonic on $M$. 

For each $h\in \mathcal{F}$, choose $a>0$ small enough and $\lambda>0$ large enough so that \eqref{hlp2} is plurisubharmonic.
Define
\begin{equation}
\label{hlp2add}
l_{h,p}(x):=\{\lambda\cdot (h^*|z|^2)(x)+h^*\big(\rho_a\ell_{P}\big)(x)\}, \quad x\in M.
\end{equation}
The function $l_{h,p} (x)$ is clearly
plurisubharmonic.
For $\epsilon >0$, define
\begin{equation*}
\varphi_{p}(x):=\sup_{h\in \mathcal{F}} \tanh( \epsilon l_{h,p}(x)),\quad x\in M.
\end{equation*}
We claim that one can find sufficiently small $\epsilon>0$ so that the bounded smooth function $l_{h,p}$ is plurisubharmonic.
Write $l=l_{h,p}(x)$. From direct computation,
we get
\begin{equation}
\sqrt{-1}\pd\op \tanh \epsilon l
=\frac4{(e^{2\epsilon l}+1)^3}\left[e^{4\epsilon l}(\sqrt{-1}\epsilon\pd\op l-2\sqrt{-1}\epsilon^2\pd l\wedge\op l)+e^{2\epsilon l}(\sqrt{-1}\epsilon\pd\op l+2\sqrt{-1}\epsilon^2\pd l\wedge\op l)\right]
\label{tanhel}
\end{equation} 
as a current.
We can choose $\lambda=\frac C{a^2}$ 
for some constant $C>0$ and sufficiently small constant $a>0$
to make sure that $l$ is plurisubharmonic.
For this purpose,  we first check the computations for $\l_1:=\lambda |h|^2$.
Now
\begin{eqnarray*}
\sqrt{-1}\epsilon\pd\op l_1-4 \sqrt{-1}\epsilon^2\pd l_1\wedge\op l_1
&=& \epsilon\lambda(1-4\epsilon\lambda|h|^2)\sqrt{-1}dh\wedge d\overline{h} \\
&\geq& \epsilon\frac{C}{a^2}(1- 4\epsilon\frac{C}{a^2})\sqrt{-1}dh\wedge d\overline{h},
\end{eqnarray*}
which is $\geq 0$ by choosing $\epsilon< \frac{a^2}{4C}$.
It follows  from the definition in \eqref{hlp2add} that the numerator of \eqref{tanhel} is at least 
$e^{4\epsilon l}(\sqrt{-1} \epsilon(\frac12\pd\op l_1 +A(\rho_a', \rho_a^{\prime\prime})))$
for some expressions $A$ in first and second derivatives of $\rho_a$ with bounded constants.  By choosing 
$a$ sufficiently small, the above expression is 
$\geq 0$ and hence $\tanh \epsilon l$ is plurisubharmonic, which implies that $\varphi_p$ is plurisubharmonic.
%
%
%
%
%
Therefore $\varphi_p(x)$ 
is a bounded plurisubharmonic function on $M$.

It remains to show that $\varphi_p$ is an exhaustion function on $M$. 
We have $\varphi_p(M)\subset [0,1)$.
Let  $\{q_n\}\subset M$ be a divergent sequence. 
Since $d_{C,M}$ is nondegenerate and complete, $d_{C,M}$ induces the complex topology of $M$ (cf. \cite[Proposition 4.1.2]{Koba1998}). We have $d_{C,M}(p,q_n)\rightarrow \infty$ as $n\rightarrow \infty$ for any fixed $p\in M$. Apply Arzela-Ascoli Theorem to each $d_{C,M}(p,q_n)$, we obtain a sequence $\{h_n\}\subset\mathcal{F}$ such that 
\[
\ell_P(h_n(q_n))=\ell_P(0, h_n(q_n))=d_{C,M}(p,q_n)\rightarrow \infty, \quad as \quad n\rightarrow \infty.
\] 
Then for each $n$,
\[
\lambda(h_n^*|z|^2)(q_n)+h_n^*(\rho_a\ell_P)(q_n)
=\lambda |h_n(q_n)|^2+ \rho_a(h_n(q_n))\ell_P(h_n(q_n))
\leq l_p(q_n).
\]
The above inequality forces $\varphi_p(q_n)\rightarrow 1$ as $n\rightarrow \infty$. 
\end{proof}

We remark that the exhaustion function obtained can be chosen to be smooth after regularization by techniques in several complex variables, such as the one given in \cite{Rich1968}, or
\cite[Proposition 1.2]{KR1981}.

\subsection{$L^2$-theory and Bergman metrics}\label{CaraL2}
In the following, we consider complex manifolds that are equipped with K\"ahler metrics.

\begin{theorem}[Andreotti-Vesentini \cite{AV1965} and H\"ormander \cite{Horm1965}]
Let $(M,g)$ be a complete K\"ahler manifold with K\"ahler form $\omega$ and Ricci form $\Ric(\omega)$. Let $(\mathcal{L},h)\rightarrow M$ be a Hermitian holomorphic line bundle with curvature form $\Theta(\mathcal{L},h)$. Write $dV_g$ as the volume form on $M$ induced by $g$ and $\|\cdot \|$ as the natural norm induced by $h$ and $g$.
Suppose $\psi$ is a smooth function on $M$ and $c$ is a positive continuous function on $M$ such that
\[
\Theta(\mathcal{L},h)+\Ric(\omega)+\sqrt{-1}\partial\barD \psi \geq c\,\omega 
\]
everywhere on $M$. Let $f$ be a $\barD$-closed square integrable $\mathcal{L}$-valued $(0,1)$-form on $M$ such that 
$\int_M \frac{\| f\|^2}{c}e^{-\psi} dV_g<\infty$. 

Then there exists a section $u$ of $\mathcal{L}$ solving the inhomogeneous problem $\barD u=f$ with the estimate
\[
\int_M \|u\|^2e^{-\psi} dV_g\leq \int_M \frac{\| f\|^2 }{c} e^{-\psi} dV_g<\infty.
\]
Furthermore, $u$ is smooth whenever $f$ is.
\label{hor}
\end{theorem}

\begin{proposition}
\label{genM}
Let $(M,g)$ be a complete K\"ahler manifold. Suppose $M$ is Carath\'eodory hyperbolic. Then the $L^2$-holomorphic sections of the canonical bundle $K_{M}$ separate points and generate jets of $K_M$ of any orders.
\end{proposition}
\begin{proof}
\noindent
\underline{Base-point freeness and generation of jets}:\\
Let $p\in M$. Given some integer $j> 0$. We want to show that there is a $L^2$-holomorphic section $\sigma\in H^0_{L^2}(M,K_M)$ such that the jets of order $j$ for $\sigma$ is nonvanishing at $p$. Consider a small neighbourhood $U\subset M$ of $p$. Let $z=(z_1,\dots,z_n)\in U$ be the local coordinates so that $p=0$. Let $\chi$ be a smooth cutoff function supported in a relatively compact neighbourhood $U_0\subset U$ of $p$. 

Let $\varphi$ be given as in Theorem \ref{psh'}.
Define
\[
\psi(z):= \alpha \varphi(z)+\chi(z) \log |z|^{2n+k},
\]
where a constant $\alpha>0$ and an integer $k\geq 0$ are to be chosen later.
Then $\psi$ defines a smooth function on $M$. Let $\omega$ be the K\"ahler form associated to the complete K\"ahler metric $g$.
Note that on $M$, $\Theta(K_M,g)=-\Ric(\omega)$, so
\begin{align*}
(*): \quad \Theta(K_M,g)+\Ric(\omega)+\sqrt{-1}\ddbar \psi 
&=\sqrt{-1}\ddbar \psi  \\
 &= \alpha \sqrt{-1}\ddbar \varphi + (2n+k)\sqrt{-1}\ddbar (\chi\log|z|).
\label{eq*}
\end{align*}
The term $\sqrt{-1}\ddbar (\chi\log|z|)$ is positive on $U_0\subset U$, vanishes on $M-U$, and has uniformly bounded negativity in $U-U_0$. Since $\varphi$ is strictly plurisubharmonic on $M$, these exists a sufficiently large $\alpha>0$ such that
\[
\alpha \sqrt{-1}\ddbar \varphi + (2n+k)\sqrt{-1}\ddbar (\chi\log|z|)\geq\epsilon \,\omega
\]
for some $\epsilon>0$. Let $v\in \Gamma(U,K_M|_U)$ be a basis. Then $\eta:=\chi v $ defines via extension by zero a smooth section of $K_M$.
Note that $\barD \eta\equiv 0$ on $U_0\cup (M-U)$. So $\supp(\barD \eta)\subset U-U_0$ is compact. Since $e^{-\psi}|_{U-U_0}$ is smooth,
\[
\int_M \frac{\|\barD \eta \|^2 }{\epsilon} e^{-\psi} dV_g<\infty.
\]
By Theorem \ref{hor}, there exists a solution $u$ to the equation 
$
\barD u= \barD \eta 
$
such that 
\begin{equation}
\label{estu}
\int_M \|u\|^2e^{-\psi} dV_g
 \leq \int_M \frac{\|\barD \eta \|^2 }{\epsilon} e^{-\psi} dV_g<\infty.
\end{equation}
Then on $U_0$, there exists a constant $C>0$ such that, 
\begin{equation}
\label{estu1}
e^{-\psi} dV_g 
= C e^{-\alpha\varphi} \cdot \frac{1}{|z|^{2n+k}}dV_e 
=C e^{-\alpha\varphi} \cdot  \frac{1}{r^{2n+k}}r^{2n-1}dS
=C e^{-\alpha\varphi} \cdot  \frac{1}{r^{k+1}}dS,
\end{equation}
 where $dV_e$ is the local Euclidean volume form, $r=|z|$ is the polar radius and $dS$ is the volume form of unit sphere. Since $\varphi$ is bounded on $M$, it follows from the estimates $\eqref{estu}$ and $\eqref{estu1}$ that $u$ vanishes at $p=0$ to the prescribed order $j\geq 0$ by choosing $k=j-1$ and $\alpha>0$ sufficiently large.

Now let $\sigma:=\eta-u$. Then $\sigma(p)=\eta(p)\neq 0$ and $\barD \sigma=0$. Hence $\sigma$ is a nontrivial holomorphic section of $K_M$ nonvanishing at $p$. Moreover, $\sigma$ is $L^2$ with respect to $e^{-\psi}=e^{-\alpha\varphi}e^{\chi\log|z|^{2n+k}}$, which is bounded from below on $M$. Thus $\sigma$ is also $L^2$ with respect to the background K\"ahler metric $g$. This shows that the $L^2$-holomorphic sections of $K_M$ is base-point free. Similarly, if we replace $\eta=\chi v$ by $\eta= z^{i_1}\dots z^{i_s}\chi v$, the above argument with minor modifications also implies that the $L^2$-holomorphic sections of $K_M$ generates jets of any given order.\\

\noindent
\underline{Separation of points}:\\
Let $p,q\in M$ be a pair of distinct points. 
To show that the $L^2$-holomorphic sections of $K_M$ separate points on $M$, 
we have to construct $\eta,\sigma \in H^0_{L^2}(M, K_M)$ such that $\eta(p)\neq 0$ and $\eta(q)=0$ while $\sigma(p)=0$ and $\sigma(q)\neq 0$. In the argument for the generation of jets, $\eta$ is obtained as desired.
If $U\subset M$ is small neighborhood of $p$ so that $q\notin U$, then we can find a neighborhood $V$ of $q$ that is disjoint from $U$.  The argument for the pair $(U,p)$ can be applied to the pair $(V,q)$ to obtain another $L^2$-holomorphic sections $\sigma$ of $K_M$ as desired. 
If $q\in U$, then take $p=0,q=z(q)\neq 0$ and replace $\psi=\alpha\varphi+\chi\log |z|^{2n+k}$ by $\psi = \alpha\varphi+\chi(\log|z|^{2n+k}+\log|z-z(q)|^{2n+k})$. The rest of the argument is then similar.
\end{proof}

\begin{lemma}
Let $(M,g)$ be a complete K\"ahler manifold. Suppose $M$ is Carath\'eodory hyperbolic. 
Then for any integer $q>0$,	the $L^2$-holomorphic sections of $qK_{M}$ separate points and generate jets of $qK_{M}$ of any orders.
\label{genqM}
\end{lemma}
\begin{proof}
Write $qK_M=(q-1)K_M + K_M$. 
Let $\{s_i\}_{i=1}^\infty $ be a unitary basis of the the $L^2$-holomorphic sections of $K_{M}$. 
The K\"ahler metric $g$ induces a metric on $K_M=\det(T^*_M)=(\det(T_M))^{-1}$, given by $(\det g)^{-1}$ or $V_g^{-1}$ in local coordinates, where $V_g$ is the coefficient of the volume form.
Define
\[
\kappa:= \left( \frac{1}{\sum_{i=1}^\infty |s_i|^2}\right )^{q-1}\cdot \det g^{-1}:=\theta\cdot \det g^{-1},
\]
which can be viewed as the product of metric $\theta$ on $(q-1)K_M$ with the metric $g$ on $K_M$. Since $\theta$ is induced by the sum of squares of holomorphic sections, we have 
$\Theta((q-1)K_M,\theta)\geq 0$
Thus the curvature $\Theta(qK_M,\kappa)\geq  \Theta(K_M,g)$, so that
\begin{align*}
(*)_q: \quad \Theta(qK_M,\kappa)+\Ric(\omega)+\sqrt{-1}\ddbar \psi 
\geq \Theta(K_M,g)+\Ric(\omega)+\sqrt{-1}\ddbar \psi .
\label{eq*q}
\end{align*}
Thus the argument in Proposition \ref{genM} for $(K_M,g)$ can be applied to $(qK_M,\kappa)$.
\end{proof}

\subsubsection*{Bergman metric}
\label{Berg}
Let $(M,g)$ be a K\"ahler manifold.
The space of holomorphic $n$-forms $\alpha$ on $M$ such that $|\int_M \alpha\wedge \overline{\alpha}|<\infty$ forms a separable complex Hilbert space $\mathcal{W}$ under the inner product $(\alpha,\beta)\mapsto \frac{i^{n^2}}{2^n}\int_M \alpha\wedge \overline{\beta}$. When $K_M$ is equipped with the natural Hermitian metric induced by $g$, $\mathcal{W}$ is equal to the space $H^0_{L^2}(M,K_M)$ of $L^2$-canonical sections with respect to $g$. The Bergman kernel on $M$ is a sum 
$
B_M:= \sum_{j=1}^\infty \frac{i^{n^2}}{2^n} s_j\wedge \overline{s_j},
$
where ${s_j}$ is a unitary basis of $H^0_{L^2}(M,K_M)$. 
The Bergman kernel $B_M$ is independent of the choice of basis.   Let 
$\tk(x):=B_M(x,x)
$ be the trace.
By abuse of notation, we also call $\tk$ the Bergman kernel on $M$.
Write $s_j=\hat{s}_jdz_1\wedge \cdots \wedge dz_n$ in local coordinates. We know that $\tk(z)=\sum_{j=1}^\infty |\hat{s}_j|^2$. 
Note that $\tk$ defines a (possibly singular) Hermitian metric on $\det T_M=K^*_M$. Write $| s_j|^2_g:=(\det g)^{-1} |\hat{s}_j|^2$ locally. Then $k(x):=\sum_{i=1}^\infty | s_i|_g^2$  defines a function on $M$, which we call the Bergman function.

The Hessian of $\log k(x)$ in local coordinates around $x$ gives rise to smooth Hermitian $(1,1)$-form that is only positive semi-definite in general. We call the corresponding pseudo-length function $\mathcal{B}_M$ the Bergman pseudo-metric. If $\mathcal{B}_M(x,v)>0$ for any $v\in T_xM-\{0\}$, then $\mathcal{B}_M$ is said to be nondegenerate at $x$. If $\mathcal{B}_M$ is nondegenerate on all of $M$, then $\mathcal{B}_M$ is called the Bergman metric.

%

It follows from Lemma \ref{genqM} that

\begin{theorem}[=Theorem \ref{bergm}]
Let $(M,g)$ be a complete K\"ahler manifold. 
Suppose $M$ is Carath\'eodory hyperbolic. Then the Bergman pseudo-metric $\mathcal{B}_{M}$ is nondegenerate. 
\label{bergm'}
\end{theorem}

Note that infinitesimal Carath\'eodory pseudo-metric decreases under holomorphic maps, i.e. if $f:Y\rightarrow X$ is a holomorphic map between complex spaces, then
$E_X(f(p),(df)_p(v))\leq E_Y(p,v)$ for any $p\in Y$ and $v\in T_pY$.
In particular,

\begin{corollary}
Let $(X,g)$ be a complete K\"ahler manifold. 
Suppose there is a holomorphic covering $\pi:X'\rightarrow X$ such that $X'$ is Carath\'eodory hyperbolic. Then for any holomorphic covering $\pi: M\rightarrow X'$, $M$ is Carath\'eodory hyperbolic and the Bergman pseudo-metric $\mathcal{B}_{M}$ is nondegenerate. 
\label{bcover}
\end{corollary}


In Corollary \ref{bcover}, since $k$ is invariant under the action of deck transformation, the Bergman metric $\mathcal{B}_{M}$ on $M$ descends to a smooth Hermitian metric $h_{X}$ on $X$. However, it is not clear if $h_{X}$ is equivalent to $\mathcal{B}_{X}$. In general, $\mathcal{B}_{X}$ may even be degenerate.
\bigskip

For any complex manifold $X$, the classical comparison theorem between Carath\'eodory pseudo-metric $E_X$ and the Bergman pseudo-metric $\mathcal{B}_X$ says that if $p\in X$ is a point such that the Bergman kernel at $p$ is nontrivial, then for any $v\in T_pX$, $E_X(x,v)\leq \mathcal{B}_X(x,v)$ (\cite{Lu1958,Hahn1976} for domains ;  \cite{Bur1977,Hahn1977,Hahn1978} for manifolds, cf. also \cite{AGK2016} ).
Under the assumption of Theorem \ref{bergm}, the nontriviality of Bergman kernel is automatic:
\begin{corollary}
Let $M$ be a complete K\"ahler manifold. Suppose $M$ is Carath\'eodory hyperbolic. Then for any $v\in T_xM$, $E_M(x,v)\leq \mathcal{B}_M(x,v)$.
\label{ineqcb}
\end{corollary}

If $M$ has many nonconstant bounded holomorphic functions but do not have a priori complete K\"ahler metrics, we may instead resort to the completeness of the Carath\'eodory distance by applying Corollary \ref{stein}:

\begin{corollary}
Let $M$ be a  strongly Carath\'eodory hyperbolic complex manifold. 
Then the Bergman pseudo-metric $\mathcal{B}_{M}$ is nondegenerate and complete. 
\label{bergnk}
\end{corollary}
It is  expected that Corollary \ref{ineqcb} and \ref{bergnk} will be useful for the study of some unbounded domains in view of \cite{AGK2016}. For example, it is well-known that for every generalized analytic polyhedron $P\subset \mathbb{C}^n$, $d_{C,P}$ is nondegenerate and $(P,d_{C,P})$ is complete, cf. \cite[Corollary 4.1.9]{Koba1998}.

%
%

%
%

\section{Hyperbolicity properties of subvarieties}

\subsection{Motivation from Kobayashi metric}
For any complex manifold $M$, one can define the Kobayashi pseudo-distance $d_{K,M}$ and the infinitesimal Kobayashi(-Royden) pseudo-metric $\mathcal{K}_M$ as in the introduction (cf. \cite{Koba1998} for detailed definitions and basic properties). In Lang's original survey \cite{Lang1986}, a complex manifold $M$ is said to be Kobayashi hyperbolic if $d_{K,M}$ is nondegenerate. Note that $d_{K,M}$ is equal to the pseudo-distance $d_{\mathcal{K}_M}$ obtained by integrating $\mathcal{K}_M$. 
So one is also lead to the study of the $\mathcal{K}_M$ when considering Conjecture \ref{langconj} and  \ref{conj}. 

In general, the infinitesimal Carath\'eodory pseudo-metric $E_M$ defines a continuous locally Lipschitz Finsler pesudo-metric with holomorphic sectional curvature $\leq -1$ \cite[p. 179-181]{Koba1998}. If $M$ is Carath\'eodory hyperbolic, then $M$ is Kobayashi hyperbolic 
and so are any manifold quotients $X=M/\Gamma$.
This shows that Carath\'eodory hyperbolicity is a stronger condition than Kobayashi hyperbolicity. 
We are going to prove that Carath\'eodory hyperbolicity is enough to imply Theorem \ref{gtcpt} and \ref{gtncpt}, which are motivated by Lang's Conjecture \ref{langconj} and its genearalization Conjecture \ref{conj} respectively.

%

\subsection{Compact manifolds}

Recall that a compact complex manifold $Z$ is said to be of general type if the Kodaira dimension $\kappa(Z)=\dim Z=:m$, or equivalently $\limsup_{q\rightarrow \infty} \frac{\dim H^0(Z,K^q_Z)}{q^m}>0$, i.e. $K_Z$ is big. In general, $\kappa(Z)\leq a(Z)\leq \dim Z=m$, where $a(Z)$ is the algebraic dimension of $Z$. A possibly singular projective algebraic variety $Z$ is said to be of general type if a smooth birational model $\tZ$ of $Z$ is of general type, (cf. \cite{Ueno1975}).

\begin{proof}[Proof of Theorem \ref{gtcpt}]
We first show that $X$ is of general type. 
In fact, we show that $K_X$ is even ample. Let $g$ be a complete K\"ahler metric on $X$.
Write $\pi: M\rightarrow X$ be the universal covering. By Theorem \ref{psh}, there exists a bounded real-analytic strictly plurisubharmonic function $\varphi$ on $M$.
Note that $\pi^*g$ is a complete K\"ahler metric on $M$.
By Theorem \ref{bergm'}
the Bergman metric is nondegenerate on $M$ and descends to a complete K\"ahler metric of strictly negative Ricci curvature on $X$. Hence by Kodaira's embedding theorem, $K_X$ is ample.  Alternatively, this follows from \cite[Section 3.2]{Gro1991} 
or the proof of \cite[Corollary 2]{Yeun2009} 
for the conclusion of general type.

Let $Z\subset X$ be a subvariety of complex dimension $m\leqslant n$.  
 Let $Y$ be an irreducible component of $\pi^{-1}(Z)$  on $M$. 
First consider the case that $Z$ is smooth.
Then the restriction of $\varphi$ to $Y$
gives rise to a bounded strictly plurisubharmonic function.
 Hence the argument of the last paragraph implies that $K_Z$ is ample and  $Z$ is of general type. 

In case $Z$ is singular, we  need to show that $\widetilde{Z}$ coming from a smooth resolution $\sigma:\tZ\rightarrow Z$ of $Z$ is of general type.  
Note that $Z$ is compact K\"ahler. The fibre product $\widetilde{Y}$ induced by the covering map $\pi|_{Y}: Y\rightarrow Z$ and the resolution $\sigma: \widetilde{Z}\rightarrow Z$ give rise to a resolution of singularity $\tau:\tY\rightarrow Y$ of $Y$, and the covering map $\pi':\tY\rightarrow \tZ$
so that $\pi\circ\tau=\sigma\circ\pi'$.

The manifold $\tY$ is K\"ahler.  
Let $Y_0$ be the smooth locus of $Y$ and $S_Y$ be the singularity of $Y$.
Note that $\tau$ restricts to give a biholomorphic map $\widetilde{Y}_0:=\widetilde{Y}-\tau^{-1}(S_Y)\cong Y-S_Y=Y_0$.
Then $\tau^*\circ (\varphi|_{Y})=\varphi\circ\tau$
is a bounded smooth plurisubharmonic function on $\tY$, which is not strictly plurisubharmonic at $\tY-\tY_0$.  The argument 
in the Proof of Proposition \ref{genM} implies that we are able to construct 
$L^2$-holomorphic sections of $K_{\widetilde{Y}}$ on $\tY$, which in fact
separate points and generate jets on $\tY_0$. 
We can now apply the technique in \cite[Section 3]{Gro1991} to construct meromorphic functions $f$ on $\tZ$  by quotients of Poincar\'e series associated to pluricanonical line bundles of $\tY$.
Recall that the Poincar\'e series associated to a $L^1$ section $F$ of $qK_{\tY}$ on $\tY$ for $q\geqslant 0$ is defined by 
\begin{equation}
\label{poin}
P_F(z)=\sum_{\gamma\in \Gamma}\gamma^*F(z),
\end{equation}
 where $\Gamma$ is the deck transformation of $\tZ$ associated to the covering map
$\pi':\tY \rightarrow \tZ$.  $P_F$ is invariant under $\Gamma$ and hence descends to give a holomorphic section of $K_{\widetilde{Z}}$ on $\tZ$.  The quotient $\frac{P_{F_1}}{P_{F_2}}$ of two such Poincar\'e series $P_{F_1}, P_{F_2}$ of the same weight $q$
gives rise to a meromorphic function $f$.

In the following, we apply the arguments from \cite[Proposition 3.2.A, Corollary 3.2.B]{Gro1991} to deduce  that the algebraic dimension of all such $f$ on $\tZ$ is 
$m$.
Since $f$ comes from quotients of Poincar\'e series which are pluricanonical sections on $\tZ$, 
for all $q>0$ sufficiently large, the image of the meromorphic map 
$\Phi_{qK_{\tZ}}$
defined by the sections of $qK_{\tZ}$ is a variety $W_q$ of dimension $\kappa(\tZ)$ (cf. \cite[Chapter 5]{Ueno1975}), and the field of meromorphic functions $\mathbb{C}(W_q)\subset \mathbb{C}(\tZ)$. It follows that there is a compact subvariety $W_q'\subset \tZ$ corresponding to $\mathbb{C}(W_q)$.
Note that there is actually no compact subvariety of dimension greater than $0$ on $Y$ as $\varphi$ is strictly plurisubharmonic.  
The only compact subvarieties on $\tY$ 
are given by the resolution of singularity of $Y$. On any relatively compact open subset $U$ of $\tY$, there are only a finite number of compact subvarieties of $\tY$ contained in $U$.
Suppose that $\kappa(\tZ)=k<m$.  Then the linear system $\Phi_{qK_{\tZ}}:\tZ\rightarrow W_q\subset \mathbb{P}^d$ has fibers $F_t, t\in W_q$, of dimension $m-k>0$, where 
$d=\dim_{\bC}(H^0(\tZ,qK_{\tZ}))$.  
Applying \cite[Corollary 3.1.A]{Gro1991}, this implies $\Phi_{qK_{\tZ}}\circ \pi'$ has generically compact fibers and hence $U$ contains infinitely many compact subvarieties of $\tY$ if $U$ is sufficiently large, contradicting the earlier observation.
We conclude that the Kodaira dimension $\kappa(\tZ)=m$ and hence $Z$ is of general type.
\end{proof}

\subsection{Noncompact manifolds}
\subsubsection{Complete noncompact K\"ahler manifolds}
In this section, we consider a noncompact complex manifold $X$.  First we prove a result which is valid for a general complete noncompact K\"ahler manifold related to the setting of the previous sections.
 \begin{theorem}
Let $X$ be a complete K\"ahler manifold.  
Assume that the universal covering $M$ of $X$ is Carath\'eodory hyperbolic. Then for any subvariety $Z\subset X$, the algebraic dimension of the space of meromorphic 
functions arising from quotients of $L^1$ holomorphic sections of the pluricanonical line bundle of $Z$ is equal to $\dim_{\bC}(Z)$.
\end{theorem} 

\begin{proof}
Let $Z\subset X$ be a subvariety. If $Z$ is compact, then by Theorem \ref{gtcpt}, we are done. In the following, we will assume $Z$ to be noncompact. 
Let $Y\subset M$ be an irreducible component of  $\pi^{-1}(Z)$ on $M$.  
In this case, the fundamental domain $\cF$ of $\tZ$ in $\tY$ is noncompact.  We define in the same way the Poincar\'e series as in \eqref{poin}
associated to $F\in H_{L^1}^0(\tY,K_{\tY}^q)$.  In the case
that $\cF$ is noncompact, we know that the infinite sum converges on compacta on $\tY$ (cf. \cite{Gro1991}).  
Hence it gives rise to a $\Gamma$-invariant holomorphic section $P_F$
of $K_{\tY}^q$.  Since 
\[
\Vert |P_F|\Vert_{\cF}=\int_{\cF}|P_F|=\int_{\cF}\sum_{\gamma\in\Gamma}|\gamma^*F|\leq \sum_{\gamma\in\Gamma}\int_{\cF}|\gamma^*F|
=\int_{\tY}|F|<\infty,
\]
we conclude that $P_F$ descends to give a $L^1$ holomorphic section of $K_{\tZ}^q$ on $\tZ$.  Again, from the work of \cite{Gro1991} as in the proof of Theorem \ref{gtcpt},
we construct a lot of meromorphic functions on $\tZ$ by taking quotients of two such Poincar\'e series.  The same proof there shows that the algebraic dimension of all such
quotients is $m$.  
\end{proof}

In order to obtain an appropriate formulation of  Lang's Conjecture \ref{langconj}  in the setting of noncompact manifold, we start with the following observation:
\begin{proposition}
There exist examples of Kobayashi hyperbolic quasi-projective manifold which contains analytic but non-quasi-projective subvarieties.
\label{egnonqp}
\end{proposition}

\begin{proof} From the results of \cite{SY1996,Siu2015}, there exists a smooth irreducible curve $C\subset P_{\bC}^2$
of high degree $d$ such that $M=P_{\bC}^2-C$ is Kobayashi hyperbolic.   For $d$ sufficiently large, the linear system $\Phi_{\cO(d)}$ associated to the powers of the hyperplane line
bundle $\cO(1)$ gives an embedding $\Phi_{\cO(d)}:P_{\bC}^2\rightarrow P_{\bC}^n$ for some $n=h^0(P_{\bC}^2,\cO(d))-1$ and mapping $C$ to a hyperplane $D$ in $P_{\bC}^n$.
Identify $M$ with $\Phi_{\cO(d)}(M)$,  which is algebraic in $\bC^n=P_{\bC}^n-D$.  Assume that $D$ is defined by $Z_{n+1}=0$ in the homogeneous
coordinates of $P_{\bC}^n$.   Let $z_i=\frac{Z_i}{Z_{n+1}}$ be the coordinates on $\bC^n=P_{\bC}^n-D$.
We can assume that the projection $\Phi_{\cO(d)}(M)\rightarrow \{z_3=\cdots=z_n=0\}\cong \bC^2$ is a finite surjective mapping so that $\Phi_{\cO(d)}(M)$ is an algebraic surface sitting over $\bC^2$.
  
Let $f:\bC^{n-1}\rightarrow \bC$ be a transcendental function.  Then the graph $G_f$ of $f$ lies in $\bC^n$ is a transcendental hypersurface in $P^n_{\bC}$, considered
as the compactification of $\bC^n$ after adding a divisor at $\infty$.  
For instance, one can let $G_f$ be given by $z_1=e^{z_2}$. Then $Y:=G_f\cap  M$ is a transcendental curve in $\Phi_{\cO(d)}(P_{\bC}^2)$, which is a compactification of
$M$.  Hence $Y$ is not a quasi-projective subvariety of $M$ while $M$ is Kobayashi hyperbolic.

\end{proof}

In view of the above examples, an appropriate formulation of Lang's Conjecture \ref{langconj} for quasi-projective manifolds should be stated in terms of log-general subvarieties and as in Conjecture \ref{conj}. Recall that a compact complex manifold $\overline{X}$ is of log-general type with respect to a divisor $D\subset \overline{X}$ if $K_{\overline{X}}+ [D]$ is big. A possibly noncompact complex manifold $X$ is said to be of log-general type if there is a compactification $\overline{X}\supset X$ such that $D:=\overline{X}-X$ is a divisor and $K_{\overline{X}}+ [D]$ is big. If $X$ is singular, then $X$ is said to be of log-general type if there is a nonsingular model $\tX$ of $X$ such that $\tX$ is of log-general type.

\subsubsection{Quasi-projective manifolds}
For the proof of Theorem \ref{gtncpt}, let us begin with some preparations.
Let $r>0$. Denote by $\Delta_r\subset \mathbb{C}$ the disk of radius $r$, $\Delta^*=\Delta-\{0\}$ the punctured unit disk and $\Delta^*_r\subset \mathbb{C}$ the punctured disk of radius $r$. Two metrics $g_1, g_2$ are said to be quasi-isometric if there exists a positive constant $c>0$ such that
$
\frac{1}{c} g_2\leq g_1\leq cg_2.
$
We denote by $g_1\sim g_2$ and $\omega_1\sim \omega_2$ for the corresponding K\"ahler forms for the quasi-isometric K\"ahler metrics. 

In the rest of the article, $X$ will always be a quasi-projective manifold. Thus we can write $X=\oX-D$ for some projective manifold $\oX$ and $D=D_1\cup\dots \cup D_s$ is a divisor of $\overline{X}$ with components in simple normal crossings.

\subsubsection*{Poincar\'e type metric and geometry of bounded type}
 
For any $x\in D\subset\overline{X}$, there is a neighborhood $\overline{U}\subset \overline{X}$ of $x$ in $\overline{X}$ with local coordinates $(z_1,\dots,z_n)$ such that $\overline{U}\cap D=\{z_1\cdots z_k=0\}$ ($1\leq k\leq n$)
 and the complement $U:=\overline{U}-\overline{U}\cap D\cong (\Delta^*)^k\times \Delta^{n-k}$.
The Poincar\'e metric $g_P$ on $U$ is defined by the associated K\"ahler form
\begin{eqnarray}
\omega_P&:=&\frac{\sqrt{-1}dz_1\wedge d\overline{z_1}}{|z_1|^2(\log|z_1|^2)^2}+\dots+\frac{\sqrt{-1}dz_k\wedge d\overline{z_k}}{|z_k|^2(\log|z_k|^2)^2}+ \sqrt{-1}dz_{k+1}\wedge d\overline{z}_{k+1}+\dots +\sqrt{-1}dz_{n}\wedge d\overline{z}_n\nonumber\\
&=&\sqrt{-1}\pd\op\psi'
,
\label{1poinform}
\end{eqnarray} 
where 
\[
\psi'=(\prod_{i=1}^k\log\log(|z_i|^{-2}))\cdot (\prod_{j=k+1}^n|z_j|^2).
\]
Let $U_{r}\cong(\Delta^*_{r})^k\times \Delta_r^{n-k}$ and $\rho$ be a smooth cutoff function supported on $U_{\frac34}$ and is identically $1$ on $U_{\frac12}$. Choose  an arbitrary smooth K\"ahler metric $\omega_0$ on $\oX$.  
Let $\omega'=a\omega_0+\sqrt{-1}\pd\op(\rho\psi')$.  Then $\omega'$ defines a complete K\"ahler metric $g'$ on $X$ for a constant $a>0$ sufficiently large.  We will also call $g'$ the Poincar\'e metric on $X$.

A  more intrinsic way to describe the Poincar\'e metric on $X$ is as follows.  Let $D_i$ be defined by $s_i=0$ locally and $h_i$ be a Hermitian metric associated to the line bundle $[D_i]$.
Then $| s_i|_h:=|s_i|^2h_i$ is a well-defined function on $X$.  Rescaling $h_i$ by a small constant if necessary, we can assume that
$| s_i|_h^2<1$ for each $i$.
Let $\psi=\prod_{i=1}^n \frac{1}{\log |s_i|_{h}^{-2}}.$  
From direct computations, on $X=\oX-D$,
\begin{eqnarray}
\sqrt{-1}\pd\op\log\psi&=&\sum_{i=1}^k\left(-\sqrt{-1}\frac{\pd\op\log | s_i|_{h}^2}{\log| s_i|_{h}^2}+
\sqrt{-1}\frac{\pd |s_i|_{h}^2\wedge\op | s_i|_{h}^2}{| s_i|_{h}^4(\log | s_i|_{h}^2)^2}\right) \nonumber \\
&=&\sum_{i=1}^k\left(\frac{\Theta([D_i],h_i)}{\log| s_i|_{h}^2}+\sqrt{-1}\frac{4  \pd | s_i|_{h}\wedge\op | s_i|_{h}}{| s_i|_{h}^2(\log |s_i|_{h}^2)^2}\right).
\label{2poinform}
\end{eqnarray}
By choosing $h_i$ sufficiently small, we can make sure that the first term of \eqref{2poinform} is bounded by $\epsilon\,\omega_0$ for any given $\epsilon>0$.  Hence 
\[
\omega:=b\omega_0+\sqrt{-1}\pd\op(\log\psi)
\] 
is positive-definite if $b$ is sufficiently large, and defines a K\"ahler metric $g$ on $X$.  It is clear from definition that
$\omega\sim\omega'$.  
We will also call $g$ the Poincar\'e metric on $X$. 
\bigskip

In the following, we always equip  $X$ with the  Poincar\'e metric $g$ defined by the K\"ahler form $\omega$.
Recall that $g$ (or $\omega$) has bounded geometry in the sense of \cite{CY80} and \cite{TY1987} (cf. \cite{Wu2008}), and satisfies the following properties:\\
(i). For any $p\in D$ and a neighbourhood $U\cong (\Delta^*)^k\times \Delta^{n-k}$ of $p$ in $X$, there exists quasi-coordinate chart given by unramified covering maps of the form
$f: W_{\frac12} =(\Delta_{\frac12})^k\times \Delta_{\frac12}^{n-k}\rightarrow U$.\\
(ii).  $(X, g)$ has finite volume.\\
(iii). All the derivatives of $g$ are bounded. Hence the  Riemannian sectional curvature of $g$ is bounded. \\

From (i),  we can cover
a neighborhood $\overline{V}$ of $D$ in $\oX$ by a finite collection of polydisks, so that the corresponding collection $\cU$ of punctured polydisks $U$ cover $V=\overline{V}-\overline{V}\cap D$ in $X$. Each $U\in \mathcal{U}$ is an image of the form $f(W_{\frac12})$.
Denote by 
$\tU=\pi^{-1}(U), \tV=\pi^{-1}(V)$ and  $\tcU=\{\pi^{-1}(U): U\in \mathcal{U}\}$ 
the pull-backs via the universal covering $\pi:M\rightarrow X$.

\begin{proposition}
\label{injrad}
Let $X$ be a quasi-projective manifold such that its universal covering $M$ has a
bounded plurisubharmonic exhaustion function. Equip $X$ with the Poincar\'e metric $g$ on $X$ given by \eqref{1poinform} or equivalently \eqref{2poinform}. Then the
injectivity radius of $g$  is bounded below by a positive constant on $M$.
\end{proposition}

\begin{proof}
Assume on the contrary that there exists a sequence $p_j\in M$ such that $\tau(p_j)\rightarrow 0$ as $j\rightarrow \infty$. Then the injectivity radius $\tau_g(\pi(p_j))\rightarrow 0$ as $j\rightarrow \infty$ on $X$. So there is a point $b\in D$ such that $\pi(p_j)\rightarrow b$ as $j\rightarrow \infty$. 

We can take $b\in D$ to be a smooth point. Let $\overline{U}\subset \overline{X}$ be an open neighbourhood of $b$ in $\overline{X}$ such that $U:=\overline{U}\cap X\cong \Delta^*\times \Delta^{n-1}$, where we have the biholomorphism $\Delta^*\cong \mathcal{D}^*\subset X$ for some disk $\Delta\cong \mathcal{D}\supset \mathcal{D}^*$ transversal to $D$ at $b$ in $\overline{X}$. We can assume $\pi(p_j)\in \mathcal{D}^*$ for all $j$.

Let $V_0\subset \pi^{-1}(U)$ be a connected component. Then $V:=\pi^{-1}(\mathcal{D}^*)\cap V_0\subset M$ is connected. Moreover,
\begin{equation}
\eta:=\pi|_{\pi^{-1}(\mathcal{D}^*)\cap V_0}: V=\pi^{-1}(\mathcal{D}^*)\cap V_0\rightarrow \mathcal{D}^*
\end{equation}
is a holomorphic covering. We obtain the universal covering 
\begin{equation}
\pi':\mathbb{H}\rightarrow V\rightarrow \mathcal{D}^*,
\end{equation}
factored as $\pi'=\pi_0\circ \eta$, where $\pi_0:\mathbb{H}\rightarrow  V$ is also a holomorphic covering.
We have the following information:
\[
\begin{tikzcd}
\mathbb{H} \arrow[rd, "\pi_0"]\arrow[dd, "\pi'"]& \\
 &  V \arrow[ld,"\eta"]\arrow[r,hook]  &M \arrow[d,"\pi"]\\
 \Delta^*\cong\mathcal{D}^*\arrow[rr, hook] &  & X
\end{tikzcd}
\]
Note that the deck transformation group of $\pi_0$ is a subgroup of $\pi_1(\mathcal{D}^*)=\mathbb Z$ and hence  as a subgroup has to be of the form $m\mathbb{Z}$ for some integer $m>0$:
\[
\begin{array}{cc}
\text{\underline{Covering map}}     &\text{\underline{Deck transformations}}\\
 \pi': \mathbb{H}\rightarrow  \mathcal{D}^*   &\mathbb{Z}: z\mapsto z+1\\ \
 \pi_0: \mathbb{H}\rightarrow  V    &m\mathbb{Z}: z\mapsto z+m\\
\eta: V\rightarrow  \mathcal{D}^* & \mathbb{Z}/m\mathbb{Z}
\end{array}
\]
There are two mutually exclusive situations, either \\
(I) $V\cong \Delta^*\cong \mathcal{D}^*$, where $m>1$; or \\(II) $V\cong \mathbb{H}$, where $m=1$.\\
%
Let $t(z)$ be the injectivity radius with respect to the Poincar\'e metric $ds^2_{\mathbb{H}}=\frac{dz\otimes d\overline{z}}{y^2}\quad (z=x+iy\in \mathbb{H})$.
In case (II), $\mathbb{H}\cong \Delta$ has injectivity radius $t(z)$ bounded from below by any real positive number as $z$ approaches the boundary, so we are done.  Hence it suffices for us to consider case (I).

Let $B\subset \Delta\cong \mathcal{D}$ be a small neighbourhood of $0$. Using the explicit map $\pi'(z)=e^{2\pi i z}$ (identifying $\mathcal{D}^*\cong \Delta^*$), we see that $(\pi')^{-1}(B\cap \mathcal{D}^*)\supset \{x+iy\in \mathbb{H}: y\geq T\}$ for some sufficiently large $T>0$. We can choose the fundamental domain $A_o\subset \mathbb{H}$ of $\Delta^*\cap B\subset \Delta^*$ for the universal covering $\pi':\mathbb{H}\rightarrow \Delta^*$  to lie in the strip $L:=\{x+iy: 0<x<1\}\subset \mathbb{H}$, i.e.,
\[
A_o=\{x+iy\in \mathbb{H}: 0<x<1, y\geq T\}.
\]
For the preimage $\eta^{-1}(B\cap \mathcal{D}^*)=\pi^{-1}(B\cap \mathcal{D}^*)\cap V_0\subset V$, the fundamental domain $A\subset \mathbb{H}$ of $\eta^{-1}(B\cap \mathcal{D}^*)\subset M$ for $\pi_0$ is a union of translations of $A_0$ by some $\gamma\in \Gamma$. 
\[
\begin{array}{ccc}
\text{\underline{Covering map}}   &\text{\underline{Subset}} &  \text{\underline{Fundamental domain}} \\
\pi': \mathbb{H}\rightarrow  \mathcal{D}^* &B\cap \mathcal{D}^* & A_o\subset L \\
\pi_0: \mathbb{H}\rightarrow  V  &\eta^{-1}(B\cap \mathcal{D}^*) & \displaystyle{A=\bigcup_{\gamma\in I\subset \Gamma} \gamma(A_0)}
\end{array}
\]

\begin{claim}
Suppose the injectivity radius of $\eta^{-1}(B\cap \mathcal{D}^*)$ is $0$, then $A$
is the union of finitely many translates of $A_0$ by $\gamma\in \Gamma$.
\end{claim}
\begin{cproof}
Note that  $\tau(p_j)\rightarrow 0$ for $p_j\in \eta^{-1}(B\cap \mathcal{D}^*) $ if and only if $t(z_j)\rightarrow 0$ as $j\rightarrow \infty$ for $z_j\in A$, where $z_j$ are liftings of $p_j$ via $\pi_0$.
Moreover,
$t(z_j)\rightarrow 0$ implies that  $z_j\rightarrow \gamma(q_\infty)
$
for some $\gamma\in \Gamma:= \pi_1(\Delta^*)\cong \mathbb{Z}$.  
We can assume $z_j\in A$ for all $j>0$.
By choosing $B$ sufficiently small, we can assume $A$ has only one cusp from the orbit $\Gamma\cdot q_\infty$.
By choosing another fundamental domain $A$ if necessary, we can assume $z_j\rightarrow q_\infty$ as $j\rightarrow \infty$. 
As a fundamental domain, $A\subset \mathbb{H}$ is connected. Thus $A$ is a tessellation of $A_0$.

Suppose $A$ is the union of infinitely many translates of $A_0=\{z=x+iy: 0<x<1, y\geq T\}$ by $\gamma\in \Gamma$. Since $\Gamma$ is generated by $z\mapsto z+1$ and $A$ is connected, $A=\{z=x+iy\in \mathbb{H}, y\geq T\}$. It follows that the covering map $\pi_0: \mathbb{H}\rightarrow V$ is a biholomorphism. 
In such case, the injectivity radius of $V$ with respect to the Poincar\'e metric is bounded below, contradicting our assumption.
\end{cproof}

It follows from the above claim that $\eta:V\rightarrow \mathcal{D}^*$ is a holomorphic covering with deck transformation group $\Gamma_0\cong \mathbb{Z}/m\mathbb{Z}$, which is a finite group.
Let $\varphi$ be a bounded plurisubharmonic exhaustion on $M$.
Let $\varphi_0(\xi):= \sum_{\gamma_0\in \Gamma_0} (\varphi|_V)(\gamma_0(\xi))$ for $\xi\in V$. Then $\varphi_0$ descends to a bounded subharmonic function on $\mathcal{D}^*$, which extends to $\mathcal{D}$. 
Note that as $\varphi_V$ is an exhaustion function, the sum $\varphi_0$ cannot be a constant.  From construction, the extended function $\varphi_0$ has maximal value at the origin.  By maximum principle for subharmonic functions, such extension must in fact be a constant function. Hence $\varphi|_V$ must be constant. However, this is a contradiction.

%
\end{proof}

\subsubsection*{Estimates of Bergman kernel}
Let $X$ be a quasi-projective manifold so that its universal covering $M$ is Carath\'eodory hyperbolic.
Equip $M$ with the K\"ahler metric obtained by pulling back the Poincar\'e metric $g$ on $X$, which is associated to the K\"ahler form $\omega$ defined by \eqref{2poinform}. By abuse of notation, the pull-back K\"ahler metric and K\"ahler form on $M$ are still denoted by $g$ and $\omega$ respectively. 
By Theorem \ref{bergm'}, the Bergman metric on $M$ is nondegenerate and descends to a smooth Hermitian metric $\omega_B$ on $X$. The Bergman kernel $\tilde{k}$ and the Bergman kernel function $k$ on $M$ also descend to $X$.

We need the technical estimates in the following proposition for our later arguments.

\begin{proposition}
Let $X$ be a quasi-projective manifold such that its universal covering $M$ supports a bounded plurisubharmonic exhaustion function.  Denote by $k$ the Bergman kernel function  on $M$.
Then there exists a constant $c>0$ such that $k(x)\leq c$ for any $x\in M$.
\label{ker-est}
\end{proposition}


\begin{proof}  
Recall that $X=\overline{X}-D$. Take $\overline{V}\supset D$ to be an open neighbourhood in $\overline{X}$ so that $V=\overline{V}-\overline{V}\cap D$. Thus $X-V$ is relatively compact. It suffices to prove the Proposition for $x\in \widetilde{V}=\pi^{-1}(V)$. Cover $V$ by open subsets of the form $U_{\frac{1}{2}}$.

By Proposition \ref{injrad},
 the injectivity radius $\tau_g$ of $g$ is positive on the universal covering $M$.    In terms of the 
coverings by quasi-coordinates, choose $0<\tau<\tau_g$ so that the geodesic ball $B_{\tau}(x)\subset M$ with respect to $g$ 
lies in a quasi-coordinate  chart of the form $W_{\frac12}$
for each $x\in \widetilde{U}_{\frac12}:=\pi^{-1}(U_{\frac12})$.
In particular, since we are considering the Poincar\'e metric $g$ on $X$, there exists $\delta>0$
so that $B_{\tau}(x)$ lies
within 
an Euclidean ball $\mathcal{U}_\delta(x)$
of radius $\delta$ in terms of the quasi-coordinates in (i). 
By choosing a constant $\tau<\tau_g$ sufficiently small,
we can assume that for each $x\in \widetilde{U}_{\frac12}$, the geodesic ball $B_{\tau}(x)$ of radius $\tau$ centred at $x$  
lies in $\widetilde{U}_{\frac12}$.

On the geodesic ball $B_\tau(x)$, it follows from our choice of quasi-coordinates in bounded geometry that there exists $c>0$ independent of $x$ such that
\begin{equation}
\frac1c<\frac{\det g(y)}{\det g(z)}<c \quad \forall y, z\in B_\tau(x).
\label{g-est}
\end{equation}
 Denote by $\|\cdot \|_g$ the norm on $H^0_{L^2}(M,K_M)$ induced by $g$.
Since Bergman kernel is independent of the choice of the unitary basis, it is well-known that we can represent the Bergman kernel by extremal functions in the sense that
\begin{equation}
k(x)=
\sup_{s\in H^0_{L^2}(M,K_M), \Vert s\Vert_g=1}
| s|_g^2= |s_x|_g^2 
\label{berg-est}
\end{equation}
for some  $ s_x\in H^0_{L^2}(M,K_M)$ with $\|s_x\|_g=1$.  Write 
$s_x(y)=\hs_x(y)dz_1\wedge\cdots\wedge dz_n$ in terms of the local coordinates on a coordinate chart and $d\mu$ the corresponding Lebesgue measure.  It follows from Cauchy's estimates  or Sub-Mean-Value Inequality that the 
Euclidean norm
\begin{equation*}
|\hs_x(x)|^2\leq \frac1{\pi \delta^2}\int_{\mathcal{U}_\delta(x)}|\hs_x(y)|^2 d\mu(y) .
\end{equation*}
Together with \eqref{g-est}, this implies that for a constant $c'>0$, 
\begin{align}
\begin{split}
|s_x(x)|_g^2&\leq  c'  \int_{\mathcal{U}_\delta(x)}|s_x(y)|_g^2 d\mu(y)  \\
&\leq c'\Vert s_x\Vert_g^2   \\
&=c'.
\label{est1}
\end{split}
\end{align}
  Hence from \eqref{berg-est} , $k(x)\leqslant  c'$  on $M$. 
  \end{proof}

Combing Proposition \ref{injrad} and Proposition \ref{ker-est} with Proposition \ref{exhaust}, we have

\begin{corollary}
Let $X$ be a quasi-projective manifold such that its universal covering $M$  is strongly Carath\'eodory hyperbolic.  Equip $X$ with the Poincar\'e metric $g$ on $X$ given by \eqref{1poinform} or equivalently \eqref{2poinform}. Then the
injectivity radius of $g$ on $M$ is positive and there exists a constant $c>0$ such that $k(x)<c$ for all $x\in M$.
\label{cinj}
\end{corollary}

\subsubsection*{Proof of Theorem \ref{gtncpt}}
\begin{proof}
Let us prove that $X$ is of log-general type. From assumption, 
$X=\oX-D$ for some divisor $D$.  By the resolution of singularity of Hironaka, we can assume that $D=D_1\cup\cdots\cup D_n$ is in simple normal crossing.
Equip $X$ with the Poincar\'e metric $g$ with associated K\"ahler form $\omega$ constructed from \eqref{1poinform} or equivalently \eqref{2poinform}.

First we claim that 
there exist a Hermitian metric $\kappa$ on $K_X$ 
and some constant $c>0$, so that
for $q$ sufficiently large, the dimension
\begin{equation}
h^0_{L^2}(X,K_X^q,\kappa,\omega)\geq c\, q^n.
\label{l2grow}
\end{equation}
As in the argument of Mok in \cite{Mok1986}, we want  to  prescribe the order of vanishing for sections of $H^0_{L^2}(X,K^q_X)$ at a point $x\in X$ up to order $c q$ for some positive constant $c$. If this can be done, then the estimates \eqref{l2grow} follows by noting that the coefficients of Taylor series expansion at $x$ of order up to $q$ has $q^n$ coefficients.

 Let $x\in X$. Let $z$ be a local coordinates around $x$ as used in equation (*) in the proof of Proposition \ref{genM}, and $\chi$ be a cutoff function supported in a small coordinate neighborhood of $x$ as described in equation (*).   Suppose we would like to prescribe the vanishing order up to $\ell$, from the notations used in equations (*) or (*q), it suffices to make sure that there exists a Hermitian metric $\kappa$ such that   
\begin{equation}
\Theta(qK_X,\kappa)+\Ric(\omega)+(2n+\ell)\sqrt{-1}\ddbar(\chi\log|z|)\geq \epsilon\, \omega
\label{needed}
\end{equation}
for some positive constant $\epsilon>0$ .

%

By construction, $\Theta(K_X,g)>0$ on a neighborhood $U$ of $D$ in $X$ and $\Theta(K_X,g)>-c_1\omega$ on $\oX-U$  for some constant $c_1>0$.  On the other hand, 
apply the proof of Theorem \ref{gtcpt}, we found that the inverse of Bergman kernel $\tk^{-1}$ defines a Hermitian metric $h_o$ on $K_M$ which descends to $K_X$, so that 
$\Theta(K_X,h_o)=-\Ric(X,h_o)=\omega_B =\sqrt{-1}\ddbar \log k(x)$.  
Here recall that $k(x)$ is the Bergman function on $M$, which descends to $X$.
Since $g$ is a metric of bounded geometry, we can choose a positive integer $\ell$ sufficiently large so that 
\begin{equation}
\ell\, \omega>c_1\Theta(K_X,g) \quad \mbox{and}
\quad  \ell\, \Theta(K_X,h_o)>2c_1\omega \quad \mbox{on $\oX-U$}. 
\label{kineq} 
\end{equation}
It follows that 
\begin{equation}
h:=h_o^{\ell/(1+\ell)}(\det g)^{-1/(1+\ell)}
\label{metric}
\end{equation}
is a smooth Hermitian metric of $K_X$, and  
\[
\Theta(K_X,h)=\frac1{\ell+1}(\ell\Theta(K_X,h_o)+\Theta(K_X,g)).
\]
Note that $\omega$ is equivalent to $\omega_P$ on $U$ defined by \eqref{1poinform}. Using \eqref{kineq},we have
\begin{eqnarray}
\Theta(K_X,h)
&\geq&\left\{\begin{array}{ll}
\frac1{\ell+1}(c_1\omega) & \mbox{on} \  \oX-U\\
\frac1{\ell+1}\Theta(K_X,g_P)& \mbox{on} \ U
\end{array}\right.\nonumber\\
&\geq&\frac{c_2}{\ell+1} \omega\ \ \mbox{on} \  X,
\end{eqnarray}
where $c_2=\min(c_1, \Theta(K_X,g_P)/\omega_P)$.
Hence there exist constants $c_3>0$ such that 
\begin{equation}
\Theta(K_X,h)>c_3\omega.
\label{hcur3}
\end{equation} 
Now from \eqref{hcur3}, by choosing $\kappa=h^q$,
\begin{eqnarray}
&&\Theta(qK_X,\kappa)+\Ric(\omega)+(2n+\ell)\sqrt{-1}\ddbar(\chi\log|z|)\nonumber\\
&\geq&(c_3q-1)\omega-c_4(2n+\ell)\omega 
\label{suff}
\end{eqnarray}
for some constant $c_4>0$.  Hence condition \eqref{needed} is satisfies for all $\ell<\frac{c_3q-1}{c_4}-2n.$  It follows that the vanishing order at $x$ can be prescribed up
to order $c_5q$ for $q$ large and the estimates \eqref{l2grow} follows.



\begin{claim}
For $q\geqslant 1$, $L^2$-holomorphic sections of $(K_{X}^q,h^q)$ over $(X,g)$ have pole  of order $\leq  q$ along $D=D_X:=\overline{X}-X$.
\end{claim}

\begin{cproof}
Let $\sigma \in H^0_{L^2}(X,K_{X}^q)$ that is $L^2$ with respect to $h^q$. By Proposition \ref{ker-est}, $h_o^{-1}\leq C\det g$ for some constant $C>0$. Thus $|\sigma|_{h}\geq |\sigma|_g$. Take a smooth point $x\in D_X$ and suppose $U\subset \overline{X}$ is a neighborhood of $x$ so that $U\cap X\cong\Delta^*\times \Delta^{m-1}$ and $U\cap D_X=\{w=0\}$, where $(w,z_2,\dots,z_m)$ is a local holomorphic coordinates of $U$. In $U$, write $\sigma=\sigma_0(dw\wedge dz_2\wedge \dots\wedge dz_n)^q$. 
%
%
Then
\begin{eqnarray*}
\infty&>&\int_{U\cap X} |\sigma|^2_{h} dV_g \nonumber\\
&\geq&\int_{U\cap X} |\sigma|^2_{g} dV_g \nonumber\\
&= &C'\cdot \int_{\Delta^*\times \Delta^{m-1}}  |\sigma_0|^2 |w|^{2(q-1)}(\log|w|)^{2(q-1)} dV_e,
\end{eqnarray*}
where $C'>0$ is a constant and $V_e$ is the Euclidean volume form.
\end{cproof}

\textit{Continuation of the Proof of Theorem \ref{gtncpt}:} It follows from the above Claim that sections of $H^0_{L^2}(X,K_X^q)$ 
that is $L^2$ with respect to $h$ obtained by solving $\barD$-equations using Theorem \ref{hor} and \eqref{suff}
extend as sections of $qK_{\oX}+qD$ over $\oX$.  
From \eqref{l2grow}, we conclude that 
\[
h^0(\oX,q(K_{\oX}+D))\geq c \, q^n
\]
for $q$ sufficiently large.
Hence $X$ is of log-general type.

Consider now $Z\subset X$ a quasi-projective subvariety.   Let $\oZ$ be the Zariski closure of $Z$ in $\overline{X}$.   If $Z$ is compact, it is a subvariety of general type from Theorem \ref{gtcpt}.
Hence we assume from this point that $\oZ\cap D\neq\emptyset.$  Let $D_Z=D\cap \oZ$.  From the result of Hironaka again, we can assume that $D_Z$ is a divisor in
simple normal crossing in $\oZ$, $D_Z=D_1\cup\cdots \cup D_k$ for some positive integer $k$, and each $D_i$ is smooth.  If $Z$ is smooth, the argument above for $X$ immediately
implies that $Z$ is of log-general type, since the restriction of the bounded strictly plurisubharmonic function $\varphi$ on $M$ as given in Theorem \ref{psh} to an irreducible component $Y$
of $\pi^{-1}(Z)$  on $M$ is still a bounded strictly plurisubharmonic function on $Y$.

Assume now that $Z$ has singularities $S_Z:=Z-Z_0$, where $Z_0$ is the set of smooth points on $Z$.  Similarly, let $S_Y:=Y-Y_0$.
Let $\sigma:\tZ\rightarrow Z$ be a smooth resolution of $Z$
coming from the resolution $\toZ$ of $\oZ$.  
 We will show that $\tZ$ is of general type.  We adopt the setting in the proof of Theorem \ref{gtcpt}.
 Let $\tau:\tY\rightarrow Y$ be the resolution of 
singularity of $Y$.   There is a covering map $\pi':\tY\rightarrow \tZ$ with $\pi\circ\tau=\sigma\circ\pi'$.
The function $\tau$ restricts to give a biholomorphic map $\widetilde{Y}_0:=\widetilde{Y}-\tau^{-1}(S_Y)\cong Y-S_Y=Y_0$, where $Y_0$ is the smooth locus of $Y$. 
Let $\varphi$ be a bounded strictly plurisubharmonic function on $M$ as given in Theorem \ref{psh}.
Then 
$\varphi\circ\tau $
 is a bounded plurisubharmonic function on $\tY$, not strictly plurisubharmonic at points of $\tY-\tY_0$.  

The key fact is the following.  Clearly
$\varphi\circ\tau $
 is strictly plurisubharmonic at a point $x\in \tY_o$.  
Since $\varphi$ is strictly plurisubharmonic in transverse direction to $\tY_o$ in $\tY$ in the sense that 
\begin{equation}
\sqrt{-1}\pd_V\pd_{\oV}\varphi(y)\geqslant c_{y,V}g(V,\oV)
\label{psh-sing}
\end{equation}
for some t $c_{y,V}>0$ depending on $y\in \tY_o$ and $V\neq0$ in transversal direction to $S_Y$.  Furthermore, if $V$ is a tangent vector of $S_Y$ at a smooth point of  $S_Y$, equation \eqref{psh-sing} holds as well.
This is reflected in the desingularization $\tY$ that for $\ty\in \tS_Y=\tau^{-1}(S_Y)$ and $0\neq W\in T_{\ty}\tY$,
\begin{equation*}
\sqrt{-1}\pd_W\pd_{\oW}\tau^*\varphi(x)\geqslant c_{x,W}\tau^*g(W,\oW),
\end{equation*}
where $c_{x,W}\geq0$ and can only be $0$ if $\tau_*(W)=0$.

The set $S_Z$ may not be codimension $1$ subvariety of $Z$.  Let $E_1$ be a divisor on $Z$ containing $S_Z$. Then $E_1$ defines a line bundle $L_1$, which has a section $s_1$ whose zero divisor $E_1$ is effective.  After replacing $Z$ by some resolution if necessary, we can assume
that $E_1$ is a union of smooth divisors in simple normal crossings.  $E_1$ may not be ample.  
Since $\oZ$ is projective, there exist a very ample line bundle $L_2$ on $\oZ$ so that it supports a holomorphic
section $s_2$ with zero divisor $E_2$, which is in simple normal crossings with $D_Z$, and that  $L_1+L_2$ is ample on $\oZ$.  Hence we can find smooth Hermitian
metrics $h_i$ on $L_i$ for $i=1,2$ such that $\Theta(L_2,h_2), \Theta(L_1,h_1)+\Theta(L_2,h_2)$ are both positive definite.
Consider the weight function $\psi$ on $\tY$ given by 
\begin{equation*}
\psi=\tau^{*}\varphi-m\tau^*(\log\Vert \pi^*s_1\Vert_{h_1}+\log\Vert \pi^*s_2\Vert_{h_2}).
\end{equation*}
It follows that 
\begin{equation*}
\sqrt{-1}\pd\op \psi=\sqrt{-1}\pd\op(\tau^{*}\varphi)+m\tau^*(\Theta(L_1,h_1)+\Theta(L_2,h_2))
\end{equation*}
is positive definite on $\tY-\tau^*(\pi^*E_1)\cup \tau^*(\pi^*E_2)$. Write $E=\tau^*(\pi^*E_1)\cup \tau^*(\pi^*E_2)$.

Applying arguments using $L^2$-estimates in \S\ref{CaraL2} with $\varphi$ replaced by $\psi$ and $M$ replaced by $\widetilde{Y}-E$,
we know that the sections
of $H^0_{L^2}(\widetilde{Y}-E, K_{\tY})$ are very ample and generate any order of jets.  Note that a section $s\in H^0_{L^2}(\widetilde{Y}-E, qK_{\tY})$
for any integer $q>0$ has to vanish on $E$, 
since 
\[
\infty>
\int_{\tY-\tau^*(\pi^*E_1)\cup \tau^*(\pi^*E_2)}|s|^2e^{\tau^*\varphi}\cdot\Vert \tau^*s_1\Vert_{h_1}^{-mq}\cdot\Vert \tau^*s_2\Vert_{h_2}^{-mq}.\
\]
It follows that $s$ extends naturally across $E$ by taking $0$ along $E$.  These sections descend along the covering $\pi':\widetilde{Y}\rightarrow \widetilde{Z}$ to sections denoted by the same notation $s$
of $H^0_{L^2}(\widetilde{Z}, qK_{\widetilde{Z}})$.  Hence the Bergman kernel $k(x)$ provides a Hermitian metric on $\tY$ and satisfies Proposition \ref{ker-est} for points on
$\tY-E$.  

This allows us to apply the earlier argument for $X$ to conclude that $h^0_{L^2}(\tZ,\tau^*K^q_{Z})\geqslant cq^m$ for $q$ sufficiently large, making use of the argument of Mok in \cite{Mok1986}.
The necessary modifications go as follows: first choose a point $x\in \tZ-\sigma^*E_1\cup\sigma^*E_2$ and metric as given in \eqref{metric}.  We need the analogue of \eqref{needed} with $X$ replaced by 
 $\tZ$ and $n$ replaced by $m$. Then choose cutoff function $\chi$ so that its support lies in $\tZ-\sigma^*E_1\cup\sigma^*E_2$.  
Use of $L^2$-estimates  as explained earlier for $X$ allows us to show that the order of vanishing at $x$ among
sections of $H^0_{L^2}(\tZ,qK_{\tZ})$ can be prescribed up to order $c_8q$ for some constant $c_8>0$.  Hence $h^0_{L^2}(Z,qK_Z)\geq c_9q^m$ 
for  $q$ sufficiently large and some constant $c_9>0$.  We conclude that $\tZ$ and hence 
$Z$ is of log-general type.
\end{proof}

\bigskip
\noindent 
Kwok-Kin Wong.
Department of Mathematics, The University of Hong Kong, Pokfulam, Hong Kong.
Email: kkwong@maths.hku.hk  \\     

\noindent
Sai-Kee Yeung.
Department of Mathematics, Purdue University, 150 N. University Street, West Lafayette, IN 47907-1395, USA.
Email: yeungs@purdue.edu

\end{document}